\documentclass[11pt]{article}
\usepackage{graphicx}
\usepackage{amssymb,amsfonts,amsthm}
\usepackage{amsmath,amsthm,amssymb}
\usepackage{amsfonts}

\newtheorem{theorem}{Theorem}
\newtheorem{lemma}[theorem]{Lemma}

\newtheorem{remark}[theorem]{Remark}
\newtheorem{definition}[theorem]{Definition}

\newcommand{\md}{{\hbox{mod }}}

\newcommand{\Cayley}{\mathrm{Cayley}}

\begin{document}

\title{On $k$-free-like groups}
\author{A.Yu. Olshanskii and M. V. Sapir
\thanks{Both authors
were supported in part by the NSF grants DMS 0455881 and
DMS-0700811.  In addition, the research of the first author was
supported in part by the Russian Fund for Basic Research grant
08-01-00573, and the research of the second author was supported
by a BSF (USA-Israeli) grant.}}

\maketitle

Let $G$ be a finitely generated group. For a subset $A \subset G$ and a finite and symmetric generating set $X$ of $G$ (i.e. $X=X^{-1}$), the set $\partial_X A = \{a\in A\colon ax \notin A$ for some $x \in X\}$ is the (inner) boundary of $A$ (relative to $X$). The Cheeger constant of $G$ with respect to $X$ is the number $\text{Che}_X G = \inf_A \frac{\# \partial_X A}{\# A}$ where $A$ runs over all nonempty finite subsets of $G$ and $\#$ denotes cardinality.

Let $k\ge 2$ be an integer. We call a group $G$ {\em $k$-free-like} if there exists a sequence of generating sets $Z_i$, $i\ge 1$, each with $k$ elements, such that the Cayley graph $\Cayley(G,Z_i)$ has girth (that is the minimal  length of a simple loop in the graph) at least $i$, and the Cheeger constant of this graph is uniformly (in $i$) bounded away from 0.

I. Benjamini (personal communication) asked whether there exists a $k$-free-like but not free group. A positive answer for $k\ge 4$ can be deduced from the paper by Akhmedov \cite{A2} and \cite{AB}:
the proof of \cite[Theorem 2.6]{A2} and the fact that non-elementary hyperbolic groups are uniformly non-amenable \cite{AB} imply that every $m$-generated non-elementary hyperbolic group is $(m+2)$-free-like. This does not provide $k$-free-like but non-free groups with $k=2,3$.

In this note, we prove the following three theorems. Theorem \ref{t1} gives, in a sense, the simplest example of $k$-free-like but not free group (for all $k\ge 2$). The proof is self-contained and relies only on the standard small cancelation theory. Theorem \ref{t0} improves \cite[Theorem 2.6]{A2} and shows that any non-elementary $m$-generated hyperbolic groups \cite{Gr} is $k$-free-like for all $k\ge m+1$, and if, in addition, the group is torsion-free, then one can take any $k\ge m$.  Theorem \ref{t} gives many $k$-free-like torsion groups. The proofs of these theorems rely on results of \cite{O93}.

\begin{theorem} \label{t1} There exist both finitely presented and non-finitely
presented 2-generated non-free groups which are $k$-free-like for
every $k\ge 2$.
\end{theorem}

\begin{theorem}\label{t0} Every non-virtually cyclic (resp. non-cyclic torsion-free) hyperbolic
$m$-generated group is $k$-free-like for every $k\ge m+1$ (resp. $k\ge m$).
\end{theorem}

\begin{theorem}\label{t}  There exists a 2-generated torsion  group
$G$ which is $k$-free-like for every $k \ge 3.$
\end{theorem}

It is obvious that a group satisfying a non-trivial identity is
not $k$-free-like for any $k$. We give a necessary and sufficient
conditions for a group to be $k$-free-like in terms of the
so-called {\em almost identities} and show that a group with
bounded girth (for all finite generating sets) may not satisfy a
non-trivial identity. The counter-examples given earlier in
\cite{Sch} and \cite{A1} were identical, and both authors
suggested to use the methods from \cite{Olbook} for the proofs;
Schleimer just ``suspected" that the approach from $\S 34$
\cite{Olbook} should work while Akhmedov claimed that $\S 34$
\cite{Olbook} explained the example. But in fact, the proof cannot
be deduced from Section $\S 34$ \cite{Olbook}, and we note at the
end of this paper that the desired examples were already provided
by Theorem 39.4 \cite{Olbook} which is not based on the techniques
from $\S 34$\cite{Olbook}.

By a result of Benjamini, Nachmias  and Peres \cite{B}, Theorems \ref{t1}, \ref{t0}, \ref{t} have applications to the theory of percolation on transitive graphs.

Let $\Gamma=(V,E)$ be a graph with vertex set $V$ and edge set $E$. An edge of the graph is called a  {\em bond}. Pick a real number $p$ between $0$ and $1$. The {\em Bernoulli bond
percolation} on $\mathcal{G}$ is a product probability measure
$\mathrm{P}_p$ on the space $\Omega=\{0,1\}^E$, the set of subsets of the
edge set $E$. For any {\em realization} $\omega \in \Omega$, the
bond $e\in E$ is called {\em open} if $\omega(e)=1$ and {\em closed}
otherwise. For $0\leq p\leq 1$ the product measure is defined via
$\mathrm{P}_p(\omega(e)=1)=p$ for all $e\in E$. Thus each bond is
open with probability $p$ independently of all other bonds.

For any realization $\omega$, open edges form a random subgraph of
$\mathcal{G}$.

The critical percolation constant  $p_c$ is defined as the infimum of all numbers $p$ for which the random subgraph has infinite connected components $\mathrm{P}_p$-almost surely.

The constant $p_c$  is one of the most important characteristics of the graph $\Gamma$, and the study of percolation on graphs, in particular Cayley graphs, has been very intensive for the last 50 years. We refer the reader to the nice article  by I. Benjamini and O. Schramm \cite{BS} and the book by R. Lyons and Y. Peres \cite{LP}.

It is known \cite{LP, BS} that the critical percolation constant of a $d$-regular tree is $\frac1{d-1}$. Hence the $p_c$ of the Cayley graph of a free group of rank $k$ with respect to its free generators is $\frac{1}{2k-1}$. It is also known \cite{BS} that the $p_c$ of a factor-graph cannot be smaller than the $p_c$ of the original graph. In particular, the $p_c$ of any $k$-generated group cannot be smaller than $\frac{1}{2k-1}$.

It is not difficult to show using \cite{LP} that the $p_c$ of a Cayley graph of a group with  respect to a generating set with $k$ generators is equal to $\frac1{2k-1}$ only if the group is free and the generators are free generators of that free group. Nevertheless, it is proved in \cite{B} that {\em a Cayley graph of a $k$-free-like group with respect to an appropriately chosen $k$-element generating set can have  $p_c$ arbitrary close to the minimal possible value: $\frac1{2k-1}$.}
Theorems \ref{t1}, \ref{t0}, \ref{t} provide plenty of examples of such groups. Notice that we do not provide any procedures of finding $p_c$ of particular Cayley graphs. Currently the only groups for which formulas for $p_c$ are known for Cayley graphs of these groups corresponding to all possible (finite) generating sets, are groups acting on simplicial trees with finite vertex stabilizers (Koz\'akov\' a \cite{Ko}). All these constants $p_c$ turned out to be algebraic numbers. All groups considered in \cite{Ko} are virtually free, hence hyperbolic. It is not known whether the $p_c$ of a Cayley graph of a hyperbolic group (relative to a finite generating set) is always an algebraic number.

For groups discussed in Theorems \ref{t1}, \ref{t}, we do not know even the cardinality of the set of possible $p_c$. It is clear from the proofs of Theorems \ref{t1}, \ref{t}, that one can construct continuously many $k$-free-like groups (for every $k\ge 2$). But it is not clear why these groups may not have the same set of critical percolation constants $p_c$.

\proof[Proof of Theorem \ref{t1}.] Recall that a symmetric  set $R$ of
cyclically reduced words which is closed under taking cyclic shifts satisfies the small cancellation
condition $C'(\lambda)$ ($\lambda
>0$) if for every two words $r,r'\in R$, $r\ne r'$,  having a
common prefix $u$, we have $|u|<\lambda\min(|r|, |r'|)$ (here
$|w|$ denotes the length of the word $w$).

We consider a finite or infinite set of positive words $R$ in
letters $a$ and $b$ which  (a) is closed under cyclic shifts, (b)
$R\cup R^{-1}$ satisfies $C'(\frac16), $ (c) $R$ has no words with
prefixes $a^2,$ $(ab)^2$ or $(ba)^2$, and (d) $R$ has no words of
length $<6.$ (Take, for example, all cyclic shifts of the words
$ab^{2j}ab^{4j}\dots ab^{100j},$ $j=1,2,\dots$) We define the
group $G = \langle a, b \mid R\rangle$.

For a given $k\ge 2,$ the set of words $X_n=X_n(k)= \{x_1=a, x_2=
ba^n, x_3=ba^{2n},\dots x_k=ba^{(k-1)n} \}$ generates $G$. We set
$n \equiv 2\ (\md 4)$.  In order to estimate the girth of $G$ with
respect to $X_n,$ consider a non-trivial cyclically reduced word
$u$ of length $\le n$ such that $u(x_1,...,x_k) =1$ in $G.$ Then
clearly $u(x_1,\dots,x_n)\ne 1$ in the free group $F(a,b)$, and so the reduced
form $U$ in generators $a,b$ of the left-hand side is non-empty.

Since $U=1$ in $G$, by the Greendlinger lemma for small
cancellation presentations \cite{LS}, $U$ must contain a subword
$V$ which is at the same time a prefix  of some $r\in R\cup
R^{-1}$ and $|V|>|r|/2\ge 3.$ Since the cyclic shifts  of $r$ have
no subwords $a^{\pm 2},$ it follows that the product
$u(x_1,\dots,x_k)$ involves some $x_j^{\pm 1}$ for $j>1.$

The factor $x_1^{\pm 1}$ occurs at most $n-1$ times in $u,$ and so
the word $U$ must have a reduced form $a^{s_0}b^{\pm
1}a^{s_1}b^{\pm 1}\dots a^{s_{t-1}}b^{\pm 1}a^{s_t},$ where
 $t\ge 1$ and $|s_1|,\dots |s_{t-1}|>0$. Since $|V|>3$ and $r$ is positive,
  the subword $V$ of $U$ must contain one of the subwords $a^{\pm 2}, (ab)^{\pm 2},
(ba)^{\pm 2} $ contrary to the choice of  $R.$ Hence the girth of
$G$ with respect to $X_n$ is at least $n+1.$

If as above, we consider a nontrivial relation of the form
$U=u(x_1^4, x_2)=1,$ then the reduced
form $U$ and any of its subwords $V$ have no occurrences of $b^{\pm 2}$ since
$n\ne 0 \ (\md 4),$ and all the exponents $s_j$ of the letter $a$ in $U$ must be even,
contrary the condition (c) for the choice of $R.$ Hence the
subgroup $\langle x_1^4, x_2\rangle$ is free in $G.$ Therefore for every
$n$ there exists a free subgroup of $G$ generated by two words of
uniformly bounded length (four) with respect to $X_n$. This implies that
the Cayley graphs $\Cayley(G,X_n)$ have Cheeger constants bounded
away from 0 (see, for example, \cite[Section 10]{AB}).

Let $R_0$ be the set of relators from $R$ with pairwise distinct
sets of cyclic shifts. It remains to note that $R_0$ is the set of
independent relators. Indeed, if one of the relations $r=1$ ($r
\in R_0$) follows from the others, then by Greendlinger lemma, $r$
must contain ``at least half" of a cyclic shift of another relator
from $R_0$ which contradicts $C'(\frac16).$\endproof

\proof[Proof of Theorem \ref{t0}.] Recall that a group $E$ is
elementary if it has a cyclic subgroup of finite index. We will
use some properties of hyperbolic groups which can be found in
\cite{O93}.

  Every hyperbolic group $H$ has a unique maximal finite normal subgroup
denoted by $E(H)$.
 If $H$ is non-elementary hyperbolic group then the quotient $H/E(H)$ is
also non-elementary hyperbolic
group.

  Let $\{a_1,\dots,a_m\}$  be a set of generators for a non-elementary
hyperbolic group $H.$
  Then one can choose $m$ pairwise distinct modulo $E(H)$ generators.

  Indeed, assume that $ a_i = a_jb $ for some $b\in E(G)$ and $i\ne j$. Then
the (images of the)
elements  $a_1,\dots,a_{j-1},a_{j+1},\dots,a_m$ generate the
infinite group $H/E(H)$. Therefore there is a word $w$ in these
$m-1$ generators such that  $w$  is not equal to any of $1,
a_1,\dots,a_m$ modulo $E(H)$, and the set $a_1,\dots a_{j-1}, wb,
a_{j+1},\dots,a_m$ generates  $H$.
 A repeated application of such a change of generators provides us with
pairwise
 distinct
generators modulo $E(H).$ Similarly, we may assume that no
generator $a_i$ belongs to the subgroup $E(H).$

  Every element $g$ having infinite order in a hyperbolic group $H$ is
contained in a unique
maximal elementary subgroup $E(g)$ \cite{Gr}. For a non-elementary
hyperbolic group $H$, (simplified versions of) Lemmas 3.4 and 3.8
from  \cite{O93} provide us with an infinite set $g_1, g_2,...$ of
elements of infinite order such that pairwise intersections of the
cyclic subgroups $\langle g_i\rangle$ generated by these elements are trivial,
 and $E(g_i) = \langle g_i\rangle \times E(H).$

  Since the subgroup $E(H)$ is finite, it follows that the cyclic
  subgroups $\langle g_i\rangle$ have pairwise trivial intersection modulo
  $E(g)$ as well.
  By the above choice
of the generators $a_1,\dots,a_m$, the set $Y=\{a_i; i\le m\}\cup
\{a_i^{-1}a_j; 1\le i<j\le m\}$ has empty intersection with
$E(H).$ Hence one can select $g=g_i$ such that the subgroup
$E(g)=\langle g\rangle  \times E(H)$ has no elements from the set $Y.$

 Given $n\ge 1$ let $X_n$ be the following generating set of $H$:
 $$X_n=\{g, a_1g^n, a_2g^{2n},\dots, a_m g^{mn}\}.$$ It is known \cite{AB} that every non-elementary hyperbolic group is uniformly non-amenable (that is the Cheeger constant is bounded away from 0 uniformly for all finite generating sets of the group). Hence in order to show that $H$ is $k$-free-like, it is enough to show that
for any given $l\ge 1$, the girth of the group $H$ with respect to
the generator set $X_n$ is greater that $l$ provided $n=n(l)$ is
large enough.

Denote by  $|x|$ the length of an element $x$ of $H$ with respect to the
generators $a_1,\dots, a_m.$
A simplified version of Lemma 2.4 \cite{O93} says that there are
$\lambda\in (0, 1]$ and $c\ge 0$ depending on $g$ only, such that
if a product of the form

 (*)   $h=x_0g^{m_1}x_1g^{m_2}x_2 \dots g^{m_t}x_t$, ($t\le l$)

 satisfies the conditions

(1) $|x_i| \le 2$ for $i=0,1,\dots,t$,

(2) $|m_i| \ge C,$ where the constant $C>0$ depends on $g$ only
($i=1,\dots,t$),

(3) $x_i^{-1}gx_i \notin E(g)$ for $i=1,\dots t-1,$

\noindent then $|h|\ge \lambda (|m_1|+...+|m_t|) - c.$

Now let us choose $n> C +\lambda^{-1}c+l$ and assume that there is
a non-trivial cyclically reduced relation of length $\le l$
between the elements of $X_n$. This relation involves at least one
generator of the form $a_ig^{ni}$ since the generator $g$ has
infinite order.

 This relation gives us an equation $h=1$
in $G,$ where $h$ has the form (*) with $x_i$ belonging to the set
$Y$, whence $h$ satisfies condition (1). It is clear from the
definition of $X_n$ that $t\ge 1$ and the condition (2) holds as
well since the generator $g$ occurs in our relation at most $l$
times and $n-l
>C.$ Finally, since every $x_i$ with $i=1,\dots, t-1$ belongs to
$Y$, we have $x_i \notin E(g)$ by the choice of $g.$ But then
$x_i^{-1}gx_i\notin E(g)$ by Lemmas 1.16 and 1.17 \cite{O93}, and the
condition (3) holds.

Thus, $|h| \ge \lambda (|m_1|+...+|m_t|) - c \ge \lambda
(n-l)-c\ge \lambda C >0$ contrary to the assumption that $h=1$ in
$G.$  The required lower bound for the girth is obtained.

If the group $H$ is, in addition, torsion free, then the
elementary subgroups in $H$ are cyclic. (It is well known that a
virtually cyclic torsion free group is cyclic.) Therefore $H$ has
a minimal generator set $\{a_1,\dots, a_m\}$ with $m\ge 2$ and
$E(a_i)=\langle a_i\rangle$ for every $i.$ In this case neither an element
$a_i$ ($1<i\le m$) nor an element $a_i^{-1}a_j$ ($1<i<j\le m$)
belongs to $E(a_1)=\langle a_1\rangle$ by the minimality in the choice of
generators. Hence the above proof works for the generator set
$S=\{a_1 (=g), a_2g^n,\dots, a_m g^{(m-1)n}\},$ and so an
$m$-generated non-cyclic torsion free hyperbolic group is
$k$-free-like for any $k\ge m.$
\endproof

To prove Theorem \ref{t}, we need the following lemma which is interesting by itself.

\begin{lemma}\label{l}
Suppose $H$ is a non-elementary hyperbolic group, and $h$ is an
element of infinite order in $H$. Then there is a natural number
$n_0=n_0(h)$ such that for every finite subset $M$ of $H$ there is
a natural number $N=N(h,M)$ such that for every $n$ divisible by
$n_0$ and $n\ge N,$ the quotient $H_1$ of $H$ over the normal
closure of $h^n$ in $H$ is non-elementary hyperbolic and the
canonical epimorphism $H\to H_1$ is injective on $M.$
\end{lemma}
\proof  The elementary subgroup $E(h)$ has a normal in $E(h)$
cyclic subgroup $C$ of finite index $n_0$. Let $C$ be generated by
an element $g$, and so $h^{n_0}=g^s$, and we may assume that
$s>0.$ By \cite[Theorem 3]{O93}, for every sufficiently large
$t=t(g,M),$ there is a canonical homomorphism $\epsilon_1$ of $H$
onto a factor-group $H_1=H/K$ such that $H_1$ is a non-elementary
hyperbolic group and $\epsilon_1$ is injective on the subset $M.$
The subgroup $K$ can be chosen as the normal closure in $H$ of the
element $g^q$ for arbitrary $q\ge q_0=q_0(g,M).$ (See line $4$ in
the proof of that theorem.) If we choose $q$ divisible by $s,$
$q=sd,$ then $g^q= h^{n},$ where $n=n_0 d.$ It now suffices to set
$N= [n_0 q_0/s]+1.$ \endproof

\proof[Proof of Theorem \ref{t}.]

Let $H=F_2 = F(a_1,a_2).$ By Osin's theorem \cite{Os2}, we can
select a large odd number $N_0$ such that the free Burnside group
$B(2,N_0)= F_2/(F_2)^{N_0}$ is uniformly non-amenable.

We can enumerate the elements of $H_0=H=\{h_1,h_2,\dots\}$ and
enumerate all pairs $(j,k),$  where $j\ge 1$ and $k\ge 3,$ such
that $i\ge j$ if $i$ is the number of a pair $(j,k).$

Assume that the canonical epimorphisms $H_0
\stackrel{\epsilon_1}{\longrightarrow} H_1
\stackrel{\epsilon_2}{\longrightarrow}  \dots
\stackrel{\epsilon_i}{\longrightarrow} H_i$ are defined for $i\ge
0,$ and (1) $H_0,\dots,H_i$ are non-elementary hyperbolic groups,
(2) (the images of) $h_1,\dots h_{i}$ have finite orders in
$H_{i}$ (3) for the pair $(j,k)$ with number $i,$ the words
 $v(i,1),...,v(i,k)$ are selected such that (a) the set
 $V_i =\{ v(i,1),...,v(i,k)\}$ generates $H_i$
 (b) $v(i,s) = a_s\ (\md (F_m)^{N_0})$  for $s =1,2$ , (c)
the girth of $H_i$ with respect to $V_i$ is at least $i$.

Consider the pair $(j,k)$ number $i+1.$  Since $H_i$ is
non-elementary hyperbolic and generated by $a_1$ and $a_2,$ we, as
in the proof of Theorem \ref{t0}, choosing an appropriate word
$g=g(i),$ can construct a generating set  $ V_{i+1} =\{
v(i+1,1)=a_1g^t, v(i+1,2)=a_2g^{2t}, \dots, v(i+1,k) =g\}$ for
$H_i$ (and for $H$ too), such that the girth of $H_i$ with respect
to $V_{i+1}$ is at least  $i+1\ge j.$ In addition we demand now
the exponent $t$ to be divisible by $N_0,$ and so $v(i+1,s) = a_s
\ (\md (F_m)^{N_0})$ for $s=1,2.$

Now we define $M$ as the set of all words in $a_1$ and $a_2$ whose
lengths are bounded from above by $d_{i+1}(i+1),$ where $d_{i+1}$
is the maximum of length of all words from $\bigcup_{l\le
i+1}V_l$. Then if (the image of)
 $h_{i+1}$ has infinite order in $H_i,$ we apply Lemma \ref{l}
 choosing the exponent $n=n(i+1)$ divisible also by $N_0.$  It provides
 us with a canonical homomorphism $\epsilon_{i+1}:H_i\to H_{i+1}$ onto
 a non-elementary hyperbolic group $H_{i+1}$ injective on the set
 $M,$ such that the kernel of $\epsilon_{i+1}$ is the normal
 closure of $h_{i+1}^n$ in $H_i.$

Such a choice of $\epsilon_{i+1}$ guarantees by induction that the
girth of $H_{i+1}$ with respect to $V_{i'}$ is at least $i'\ge j'$
if a pair $(j',k')$ has number $i'\le i+1.$

If the word $h_{i+1}$ has finite order in $H_i,$ then we set
$H_{i+1} = H_i,$ and $\epsilon_{i+1}$ is the identical mapping. In
any case the images of $h_1,\dots, h_{i+1}$ in $H_{i+1}$ are of
finite orders.  Hence, the limit group $\hat H$ for the sequence
of epimorphisms $\epsilon_i$ is a torsion group. This group is
infinite since non-elementary groups $H_i$-s are infinite.

Furthermore, the relations from every finite set of $\hat
H$-relations follow from the relations of some group $H_i,$ and so
the girth of $H$ with respect to arbitrary set $V_{i'}$ is at
least $j+1$ if $i'$ is the number of a pair $(j,k),$ since this
property holds for every $H_i$ with $i\ge i'.$

Finally, the group $\hat H$ can be canonically mapped onto
$B(2,N_0)$ by the choice of exponents $n(i)$ divisible by $N_0.$
Under this mapping, the words $v(i,1)$ and $v(i,2)$ are mapped to
the generators $a_1, a_2$ of the group $B(2,N_0)$ since the other
factors of $v(i,1), v(i,2)$ vanish. Since by \cite{Os2}, the
Cheeger constants for the images of the sets $V_i$ are uniformly
separated from $0,$ the same is true for the generating set $V_i$
of the group $\hat H$ (see \cite{AB}). Thus, the theorem is proved
for the group $G=\hat H.$
\endproof

\begin{remark}\label{r1}\rm{ Using results of Osin \cite{Os1, Os3}, one can replace ``hyperbolic" in Theorem \ref{t0} by ``strongly relatively hyperbolic". In particular, for every two finitely generated groups $A, B$ of orders $\ge 3$, the free product $A*B$ is $k$-free-like for all sufficiently large $k$. Note that if $A$ or $B$ is not finitely presented, $A*B$ is also not finitely presented.}
\end{remark}

\begin{remark}\label{r2} \rm {By choosing the exponents $n(i)$ in the proof of Theorem \ref{t} large enough, one can ensure that the group $G$ in that theorem is {\em lacunary hyperbolic}\cite{OOS}, i.e. one of its asymptotic cones is an $\mathbb{R}$-tree. More complicated but similar in spirit constructions of torsion lacunary hyperbolic groups can be found in \cite[Section 6]{OOS}.} \end{remark}

\begin{definition} A word $u$ in $k$ variables is called a $k$-almost identity of a $k$-generated group
$G$ if $u(a_1,...,a_k)=1$ in $G$ for every generating set $\{a_1,...,a_k\}$ of $G$.
\end{definition}

For example, the words  $a_1^2 a_2^2$ and $a_1^2a_2a_1^2a_2^{-1}$
are 2-almost identities but not identities of the quaternion group
of order $8$ and, respectively, of the symmetric group $S_3.$ The
left-hand sides of all $k$-almost identities of the group $G$ form
a characteristic subgroup $C=C(G)$ of the free group $F_k$ since
$C$ is the intersection of the kernels of all epimorphisms $F_k\to
G.$

Clearly if a group $G$ has a non-trivial $k$-almost identity, then the girth of $G$ with respect to every $k$-element generating set is bounded from above. The next Theorem shows that the converse statement holds too.

\begin{theorem}\label{r4} Let $G$ be a $k$-generated group. The girth of the group $G$ with respect
to every $k$-element generating set is uniformly bounded if and only if $G$ satisfies  a non-trivial
$k$-almost identity.
\end{theorem}

\proof We only need to prove the ``only if" implication. Suppose
that for some $k\ge 2$, the girth of a group $G$ with respect to
every system of generators  $\{a_1,..., a_k\}$ does not exceed
$N$, that is  $w(a_1,...,a_k)=1$ for some non-trivial in $F_k$
 word in $k$ variables of length $\le N$. Let $W=\{w_1,...,w_M\}$
be the (finite) set of non-trivial words of length $\le N$ in $k$
variables. Consider a sequence of words word $u_1,...,u_M$
constructed by induction. Let $u_1=w_1$. Suppose that we already
have $u_{i-1}$. If $u_{i-1}$ commutes with $w_i$ in the free
group, i.e. $u_{i-1}^s=w_i^t$ for some $s,t\ne 0$, we set
$u_i=u_{i-1}^s$, otherwise we set $u_i=[u_{i-1}, w_i]$. Then the
word $u=u_M$ is non-trivial in the free group but
$u(a_1,...,a_k)=1$ in $G$ for every generating set
$\{a_1,...,a_k\}$ of $G$ (and in fact for every $k$-tuple of
elements $a_1,...,a_k$ such that $w(a_1,...,a_k)=1$ for some $w\in
W$). Therefore $u=1$ is a non-trivial $k$-almost identity of $G$.
\endproof

Finally we provide an example of a group satisfying a $k$-almost
identity, but containing a free non-abelian subgroup and thus does
not satisfy any non-trivial identity.

\proof Let $n>1$ be an odd integer. In the free group $F_2=\langle a,b\rangle$,
we choose the subgroup $^nF_2$ which is generated by all $n$-th powers of words $w$ such that
the total (algebraic) number of occurrences of either $a$ or $b$ in $w$ is not divisible by $n$.
Clearly, $^nF_2$ is normal in $F_2$. Let $G=F_2/^nF_2$. \cite[Theorem 39.4]{Olbook} states that
if $n$ is large enough, the group $G$ contains a free non-abelian subgroup.
On the other hand, if words $x_1,...,x_k$ represent elements in $G$ that generate $G$,
then the number of occurrences of either $a$ or $b$ in one of $x_i$ is not divisible by $n$.
Indeed, otherwise all $x_i$ would be in the kernel of the natural
homomorphism $G\to \mathbb{Z}/n\mathbb{Z}\times \mathbb{Z}/n\mathbb{Z}$.
This $x_i$ must satisfy $x_i^n=1$ in $G$. Hence the girth of $G$ with respect to the
generating set $\{x_1,...,x_k\}$ does not exceed $n$. \endproof

\noindent Alexander Yu. Olshanskii:\\
{\small \sc Department of Mathematics, Vanderbilt University , Nashville, TN 37240.\\
Department of Mathematics, Moscow State University, Moscow, 119899, Russia.\\}
{\it E-mail:} {\tt alexander.olshanskiy@vanderbilt.edu}

\vspace{3mm}

\noindent Mark V. Sapir:
{\small\sc \\ Department of Mathematics, Vanderbilt University, Nashville, TN 37240.\\}
{\it E-mail: } {\tt m.sapir@vanderbilt.edu}


\begin{thebibliography}{1}
\label{bibbb}
\bibitem[Akh1]{A1} Azer Akhmedov,
On the girth of finitely generated groups.
J. Algebra 268 (2003), no. 1, 198--208.

\bibitem[Akh2]{A2} Azer Akhmedov,
The girth of groups satisfying Tits alternative.
J. Algebra 287 (2005), no. 2, 275--282.


\bibitem[ABLRS]{AB}G. Arzhantseva, J. Burillo,M.  Lustig, L. Reeves, H. Short, E. Ventura,
Uniform non-amenability.
Adv. Math. 197 (2005), no. 2, 499--522.

\bibitem[BNP]{B} Itai Benjamini, Assaf Nachmias, Yuval Peres, Is the critical percolation probability local? preprint, 2008.

\bibitem[BS]{BS}
Itai Benjamini and Oded Schramm, Percolation beyond {$\mathbf{Z}\sp d$},
  many questions and a few answers, Electron. Comm. Probab. \textbf{1} (1996), no.\ 8, 71--82 (electronic).

\bibitem[Gr]{Gr}
M.Gromov,
Hyperbolic groups, in: Essays in Group Theory
(S.M.Gersten, ed.), M.S.R.I. Pub. 8, Springer, 1987, 75--263.

\bibitem[Ko]{Ko} Iva Koz\'akov\'a, Critical percolation on Cayley graphs of groups acting on trees, preprint,  arXiv:0801.4153, 2008.

\bibitem[LS]{LS} Roger Lyndon and Paul Schupp.
\newblock {\it Combinatorial group theory}.
\newblock Springer-Verlag, 1977.

\bibitem[LP]{LP}
Russell Lyons and Yuval Peres, \emph{Probability on trees and networks},
  Cambridge University Press,
  http://mypage.iu.edu/$\sim$rdlyons/prbtree/prbtree.html, To appear.

\bibitem[Ol91]{Olbook}A. Yu. Olshanski, Geometry of Defining Relations in Groups, Kluwer Academic, 1991.

\bibitem[Ol]{O93}
A.Yu. Olshanskii. On residualing homomorphisms and $G$--subgroups of
hyperbolic groups. Int. J. Alg. Comp., {\bf 3} (1993), 4, 365--409.




\bibitem[OOS]{OOS} A. Yu. Olshanskii, D. V. Osin, M. V. Sapir, Lacunary hyperbolic groups, arXiv, math/0701365, 2007.

\bibitem[Os06]{Os1}
D.Osin, Relatively hyperbolic groups: Intrinsic geometry, algebraic
properties, and algorithmic problems, Memoirs Amer. Math. Soc., 179
(2006), no. 843.

\bibitem[Os07]{Os2} D. Osin, Uniform non-amenability of free Burnside groups.  Arch. Math. (Basel)  88  (2007),  no. 5, 403--412.

\bibitem[Os07']{Os3} D. Osin,
Peripheral fillings of relatively hyperbolic groups.
Invent. Math. 167 (2007), no. 2, 295--326.

\bibitem[Sch] {Sch} S. Schleimer, On the girth of groups.
Preprint.

\end{thebibliography}
\end{document}